\documentclass[11pt]{article}
\usepackage{amsmath,amsthm,amssymb,xypic}

\theoremstyle{plain}  

\newtheorem{thm}{Theorem}[section]
\newtheorem{prop}[thm]{Proposition}
\newtheorem{lem}[thm]{Lemma}

\theoremstyle{definition}

\newtheorem{prop-defn}[thm]{Proposition--Definition}
\newtheorem{defn}[thm]{Definition}
\newtheorem{notn}[thm]{Notation}

\theoremstyle{remark}

\newtheorem{rem}[thm]{Remark}

\DeclareMathOperator{\Proj}{Proj}

\DeclareMathOperator{\Bl}{Bl}

\DeclareMathOperator{\Hom}{Hom}

\DeclareMathOperator{\red}{red}

\DeclareMathOperator{\Hilb}{Hilb}
\DeclareMathOperator{\Chow}{Chow}

\DeclareMathOperator{\Lie}{Lie}

\DeclareMathOperator{\conv}{conv}

\DeclareMathOperator{\SL}{SL}
\DeclareMathOperator{\Cone}{Cone}
\DeclareMathOperator{\Aff}{Aff}
\DeclareMathOperator{\Wt}{Wt}

\newcommand{\QED}{\ifhmode\unskip\nobreak\fi\quad {\rm Q.E.D.}} 

\newcommand{\bC}{\mathbb C}

\newcommand{\bP}{\mathbb P}
\newcommand{\bQ}{\mathbb Q}
\newcommand{\bR}{\mathbb R}

\newcommand{\bZ}{\mathbb Z}

\newcommand{\cB}{\mathcal B}

\newcommand{\cX}{\mathcal X}

\newcommand{\sD}{\mathcal D}

\newcommand{\cI}{\mathcal I}
\newcommand{\cO}{\mathcal O}

\newcommand{\cN}{\mathcal N}

\newcommand{\cQ}{\mathcal Q}

\newcommand{\cS}{\mathcal S}

\newcommand{\cZ}{\mathcal Z}

\newcommand{\map}{\rightarrow}
\newcommand{\da}{\downarrow}
\newcommand{\inj}{\hookrightarrow}

\newcommand{\cq}{/ \!  /}
\newcommand{\hq}{/ \! / \! /}

\title{Compact moduli of hyperplane arrangements}
\author{Paul Hacking}
\begin{document}
\maketitle

\begin{abstract}
The minimal model program suggests a compactification of the moduli space of hyperplane arrangements which
is a moduli space of stable pairs. Here, a stable pair consists of a scheme $X$ which is a degeneration of
projective space and a divisor $D=D_1+\cdots+D_n$ on $X$ which is a limit of hyperplane arrangements. 
For example, in the 1-dimensional case, the stable pairs are stable curves of genus $0$ with $n$ marked points.
Kapranov has defined an alternative compactification using his Chow quotient construction, which may be described
fairly explicitly. We prove that these two compactifications coincide. We deduce a description of all stable pairs.
\end{abstract}

\section{Introduction}

Hyperplane arrangements have been the subject of intensive study in algebraic geometry, combinatorics and topology (see, e.g., \cite{OT}).
Applications of the theory include generalised hypergeometric functions, representations of braid groups, conformal field theory, etc.
The minimal model program provides a natural compactification of the moduli space of generic hyperplane arrangements which is a 
moduli space of certain `generalised arrangements' called stable pairs. 
We prove that this compactification coincides with an earlier compactification defined by Kapranov using a quotient construction.
Using this identification, we describe the degenerate stable pairs rather explicitly.

For us, a \emph{hyperplane arrangement} is an ordered collection $M_1,\ldots,M_n$ of $n$ hyperplanes in complex projective space 
$\bP^{k-1}$. It is in \emph{general position} if the divisor $M_1+\cdots+M_n$ is a normal crossing divisor on $\bP^{k-1}$. 
We say two arrangements are isomorphic if they are related by an automorphism of the ambient projective space, and 
denote the moduli space of arrangements in general position by $H_{k,n}$. 

We regard a hyperplane arrangement $M_1,\ldots,M_n \subset \bP^{k-1}$  as a pair  
consisting of a smooth variety $\bP^{k-1}$ and a (labelled) divisor $M_1+\cdots+M_n$. 
Then, assuming the minimal model program in dimension $k$, there is a compactification of $H_{k,n}$ which is a moduli 
space of \emph{stable pairs} \cite{A2}.
A stable pair is a pair consisting of a possibly reducible scheme $X$ and a divisor $D$ 
such that $(X,D)$ is mildly singular, and satisfies a certain `stability condition' 
which ensures that it has only finitely many automorphisms. 
In our context, there is a decomposition $D=D_1+\cdots+D_n$ and $(X,D)$ admits a smoothing
to a hyperplane arrangement.

Unfortunately, the minimal model program is only complete in dimensions $k \le 3$,
so the compactness of the moduli space of stable pairs is conjectural for $k>3$.
Moreover, the construction of the moduli space is highly abstract, and the task
of explicitly classifying the stable pairs or describing the moduli space appears
formidable.

The Chow quotient $C_{k,n}=({\bP^{k-1}}^{\vee})^n\cq \SL(k)$ provides an alternative compactification of 
$H_{k,n}$ \cite{K1}.  
For $X$ a complex projective variety and $G$ an algebraic group acting on $X$,
the \emph{Chow quotient} $X \cq G$ is the closure of the locus of generic orbit closures 
$\overline{Gx} \subset X$ in the Chow variety. It is a compactification of a geometric quotient $U/G$,
where $U$ is a Zariski open subset of $X$.
We note that Chow quotients are often better suited to moduli problems than GIT quotients. 
For example, the Chow quotient $(\bP^1)^n \cq \SL(2)$ is identified with $\overline{\mathfrak{M}}_{0,n}$
(\cite{K1}, p.~98, Theorem~4.1.8), whereas the various GIT quotients (determined by choosing a linearisation
of the action) have boundary points corresponding to point configurations $P_1,\ldots,P_n \in \bP^1$
in which some of the points coincide. In general, the Chow quotient dominates each GIT quotient, hence provides
a finer compactification. 

We prove that the Chow quotient $C_{k,n}$
is identified with the moduli space of stable pairs, at least set-theoretically.
\begin{thm}
There is a family of stable pairs over the Chow quotient $C_{k,n}$
such that the fibre over a point of $H_{k,n}$ is the corresponding hyperplane arrangement,
and no two fibres are isomorphic.
Assuming the minimal model program in dimension $k$, the family induces a bijection of sets
$$C_{k,n} \stackrel{\sim}{\longrightarrow} \{ \mbox{ stable pairs } \} / \cong$$ 
In other words, each stable pair occurs exactly once as a fibre of the family.
\end{thm}
Recall that the minimal model program is complete only in dimensions $k \le 3$, so the full result
is only proven for pairs of dimension $k-1 \le 2$. However, the MMP is only required to establish 
separatedness of moduli of stable pairs, and it may be possible to prove this by another method 
(cf. Remark~\ref{MMP}).

We deduce a fairly explicit description of all stable pairs. 
First, the intersections of the irreducible components of $X$ 
and $D$ define a stratification of the variety $X$. The poset of strata encodes the `combinatorial type' of $(X,D)$.
The possible combinatorial types correspond to polyhedral subdivisions of a certain lattice polytope. 
Moreover, each stratum is a normal, rational variety.
Second, the singularities of the pair $(X,D)$ are toric. That is, $X$ is locally obtained by glueing normal 
toric varieties along codimension $1$ toric strata, and the divisor $D$ is given by the remaining codimension $1$
strata. Moreover, a combinatorial analysis allows a complete classification of the singularities.
For example, if $k=3$ there are precisely 10 isomorphism types of germs $(P \in X,D)$.
Finally, there is a canonical embedding of $(X,D)$ in the Grassmannian of $(k-1)$-planes in $(n-1)$-space,
defined using a vector bundle of `logarithmic differentials' on $(X,D)$. Roughly speaking, the logarithmic 
differentials are 1-forms on $X$ with logarithmic poles along the divisor $D$ and the singular locus.

The Chow quotient $C_{k,n}$ is singular in general for $k \ge 3$.
For example, $C_{3,6}$ has $15$ singular points of type $(0 \in C(\bP^1 \times \bP^2))$, 
where  $C(\bP^1 \times \bP^2)$ is the cone over the Segre embedding $\bP^1 \times \bP^2 \inj \bP^5$.
However, we expect that the singularities are toric, i.e.,
the space $C_{k,n}$ together with its boundary is locally isomorphic to a toric variety together with 
its toric boundary. For, we define a canonical log structure on each stable pair such that the log deformations 
are locally unobstructed, and $C_{k,n}$ parametrises certain log deformations of stable pairs.
If the log deformations are globally unobstructed, it follows that $C_{k,n}$ parametrises all log deformations 
of stable pairs and has toric singularities as claimed. Details will appear in a subsequent paper.

We describe the case $k=3$, $n=5$ ($5$ lines in $\bP^2$).
We freely reorder the divisors $D_1,\ldots, D_5$ for simplicity of notation.
There are just two types of degenerate surfaces $X$.
The first surface $X$ has two components $X_1 \cong \bP^2$ and $X_2 \cong \Bl_P \bP^2$, and the exceptional curve 
$E \subset X_2$ is glued to a line $\Delta_1 \subset X_1$.
Let $p \colon X_2 \map \bP^1$ denote the $\bP^1$-fibration, with fibres the strict transforms of 
the lines through $P \in \bP^2$.
Each of the divisors $D_1,D_2,D_3$ is a Cartier divisor consisting of a line on $X_1$ together with 
the fibre of $p$ which meets the double curve in the same point. The divisors $D_4$ and $D_5$ are sections
of $p$ which are disjoint  from the negative section $E$, i.e., strict 
transforms of lines in $\bP^2$ not passing through $P$.
The second surface has three components $X_1,X_2 \cong \bP^2$ and $X_3 \cong \bP^1 \times \bP^1$.
Let $p_1, p_2 \colon X_3 \map \bP^1$ denote the two projections.
For $i=1$ and $2$, a line $\Delta_i \subset X_i$ is glued to a fibre  $f_i \subset X_3$ of $p_i$.
The divisor $D_1$ consists of two lines, one on each of $X_1$ and $X_2$, which both pass through
the common point $f_1 \cap f_2$ of $X_1$ and $X_2$. 
Each of the divisors $D_2,D_3$ is a Cartier divisor consisting of a line on $X_1$ together with
the fibre of $p_2$ meeting $\Delta_1$ at the same point; similiarly $D_4$ and $D_5$ are Cartier
divisors consisting of a line on $X_2$ and a fibre of $p_1$.
In each case, the condition on the singularities of $D \subset X$ may be stated as follows:
no component of the divisor $D$ is contained in the singular locus of $X$, and
$D$ has normal crossings away from the singular locus.
\\

\noindent\emph{Acknowledgements:} This paper relies heavily on earlier work of Kapranov \cite{K1}.
I would like to thank Igor Dolgachev for suggesting the problem 
of compactifying moduli of point configurations in projective space, and for
subsequent helpful conversations. Valery Alexeev has informed me that he has independently obtained some of the 
results proved here.

\section{Statement of the Main Theorem}

\begin{defn} A \emph{hyperplane arrangement} is an ordered collection $M_1,\ldots,M_n$ of hyperplanes in $\bP^{k-1}$.
We say $M_1,\ldots,M_n$ are in \emph{general position} if $M_{i_1} \cap \cdots \cap M_{i_k} = \emptyset$ for
each $i_1<\cdots<i_k$, equivalently, $M_1+\cdots+M_n$ is a normal crossing divisor on $\bP^{k-1}$.
\end{defn}

\begin{defn}
Fix $n >k$.
Let $U \inj ({\bP^{k-1}}^{\vee})^n$ denote the open locus of arrangements in general position.
Define $H_{k,n}=U/ \SL(k)$, the moduli space of arrangements in general position.
\end{defn}

We give the precise definition of a stable pair. We first define semi log canonical singularities,
the allowed singularities on a stable pair $(X,D)$. In fact, the semi log canonical singularities 
occurring on a stable pair are rather special (see \ref{slcremark}).

\begin{defn}
Let $X$ be a reduced scheme and $D$ an effective Weil divisor on $X$.
The pair $(X,D)$ is \emph{semi log canonical} if
\begin{enumerate}
\item The scheme $X$ is Cohen--Macaulay and has only normal crossing singularities in codimension $1$.
\item Let $K_X$ denote the Weil divisor class on $X$ corresponding to the dualising sheaf $\omega_X$.
Then the divisor $K_X+D$ is $\bQ$-Cartier.
\item Let $\nu \colon X^{\nu} \map X$ be the normalisation of $X$. Let $\Delta$ denote the double locus of $X$
and write $D^{\nu}$ and $\Delta^{\nu}$ for the inverse images of $D$ and $\Delta$ on $X^{\nu}$.
Then the pair $(X^{\nu},\Delta^{\nu}+D^{\nu})$ is log canonical.
\end{enumerate}
\end{defn}

\noindent We use the abbreviation slc for semi log canonical.

\begin{rem}
\begin{enumerate}
\item The dualising sheaf $\omega_X$ is $S_2$. It is also invertible in codimension $1$ by (1).
Hence it corresponds to a Weil divisor class $K_X$ as stated. If $X$ is normal, this is 
of course the usual canonical divisor class. In general $\nu^{\star}K_X = K_{X^{\nu}}+\Delta^{\nu}$.
\item No component of the divisor $D$ is contained in the double locus of $X$ by (3).
\end{enumerate}
\end{rem} 

\begin{defn} \label{def_stable_pair}
A \emph{stable pair} $(X,D)$ of type $(k,n)$ is a pair consisting of a proper scheme $X$ of 
dimension $k-1$ and a Weil divisor $D=D_1+\cdots+D_n$ on $X$ such that
\begin{enumerate}
\item The pair $(X,D)$ is slc.
\item The $\bQ$-Cartier divisor $K_X+D$ is ample.
\item The pair $(X,D)$ admits a smoothing to a hyperplane arrangement $(\bP^{k-1},M_1+\cdots+M_n)$. 
More precisely, there is a flat family of pairs $(\cX,\sD=\sD_1+\cdots+\sD_n)/T$ over the germ
of a curve such that the  special fibre is $(X,D)$, the general fibre is a hyperplane arrangement
and the divisor $K_{\cX}+\sD$ is $\bQ$-Cartier.
\end{enumerate}
\end{defn}

\begin{rem} \label{slcremark}
Assuming the minimal model program in dimension $k$, we can describe the possible singularities
of a stable pair $(X,D)$ more specifically as follows (cf. Proposition~\ref{(X,D)slc}).
For $Y \subset X$ an irreducible component, let $\Delta_Y$ and $D_Y$ denote the restriction of 
the double locus and the divisor $D$. Then the pair $(Y,\Delta_Y+D_Y)$ is locally isomorphic to 
a toric variety together with its toric boundary.
The pair $(X,D)$ is obtained by glueing the toric components $(Y,\Delta_Y+D_Y)$ 
along the codimension $1$ toric strata $\Delta_Y$. 
In particular, the divisor $K_X+D$ is Cartier.
\end{rem}

We introduce the compactification of $H_{k,n}$ obtained by the Chow quotient construction.
We briefly review the definition of the Chow quotient (cf. \cite{K1}, p. 34).
Let $X$ be a complex projective variety and $G$ a reductive algebraic group 
acting on $X$. For each $x \in X$ consider the orbit closure $\overline{Gx} \subset X$,
a closed subvariety.
There exists a Zariski open subset $U \subset X$ such that for each $x \in U$,
the subvariety $\overline{Gx}$ has fixed dimension $r$ and represents a fixed homology class 
$\delta \in H_{2r}(X, \bZ)$. We may assume $U$ is $G$-invariant.
Let $\Chow_r(X,\delta)$ denote the Chow variety of $r$-dimensional algebraic
cycles in $X$ with homology class $\delta$ \cite{B}. We have an embedding
$$ U/G \inj \Chow_r(X,\delta), \ x \mapsto [\overline{Gx}].$$
We define the \emph{Chow quotient} $X\cq G$ to be the closure of $U/G$ in $\Chow_r(X,\delta)$.
Note that it does \emph{not} depend on the choice of $U$.

In our case $X=({\bP^{k-1}}^{\vee})^n$, $G=\SL(k)$, and we may take $U$ to be the locus of 
hyperplane arrangements in general position, as above.
\begin{defn}
Define $C_{k,n}=({\bP^{k-1}}^{\vee})^n \cq \SL(k)$, the Chow quotient compactification of $H_{k,n}=U/ \SL(k)$. 
\end{defn}

\begin{thm}[Main Theorem] \label{mainthm}
There is a family of stable pairs over the Chow quotient $C_{k,n}$ such that the fibre over a point of $H_{k,n}$
is the corresponding hyperplane arrangement, and no two fibres are isomorphic.
Assuming the minimal model program in dimension $k$, the family induces a bijection of sets
$$C_{k,n} \stackrel{\sim}{\longrightarrow} \{ \mbox{ stable pairs } \} / \cong$$ 
In other words, each stable pair occurs exactly once as a fibre of the family.
\end{thm}

Morally speaking, the Chow quotient $C_{k,n}$ \emph{is} the moduli space
of stable pairs. However, the definition of an allowable family of stable
pairs over an arbitrary base is unclear, hence the scheme theoretic 
structure on the moduli space is not defined at present.
The problem is that, for $k \ge 3$, there are stable pairs that admit deformations to non-smoothable pairs.
The deformation theory in the simplest cases is analogous to that of degenerate K3 surfaces described by Friedman
(\cite{Fr}, p.~108, Theorem~5.10). These deformations must be ruled out in order 
to obtain a well defined scheme structure on the moduli space of (smoothable) stable pairs. 
We hope to resolve these difficulties using the canonical log structure on each stable pair 
defined in Section~\ref{logstr} (cf. \cite{KN}).

\section{Construction of the universal family of stable pairs} 

We construct the family of stable pairs over the Chow quotient $C_{k,n}$.
We first identify $C_{k,n}=({\bP^{k-1}}^{\vee})^n \cq \SL(k)$ with the Chow quotient $G(k,n)\cq H$
of a Grassmannian by a torus via the Gel'fand-McPherson isomorphism. 
We give a careful analysis of the universal family of cycles over $G(k,n)\cq H$
refining earlier results of Kapranov (\cite{K1}, Section~1.2) and Lafforgue \cite{L}.
We use a toric description of the universal family over the Chow quotient of a projective
space by a torus, essentially due to Lafforgue, which may be of independent interest.
Finally, we obtain the family of stable pairs as a transversal slice of the universal family
of cycles over $G(k,n)\cq H$. We deduce a description of the singularities and combinatorial types of
stable pairs.

\subsection{The Gel'fand--McPherson isomorphism}

Let $U \subset ({\bP^{k-1}}^{\vee})^n$ denote the locus of hyperplane arrangements
in general position as above.
Let $G(k,n)$ denote the Grassmannian of $k$-planes in $\bC^n$ and $G^0(k,n) \inj G(k,n)$
the locus of $k$-planes $L$ such that each Plucker coordinate $P_{I}$ (where $|I|=k$, $I \subset [n]$) is 
nonzero. Let $H={\bC^{\times}}^n/\bC^{\times}$ be the quotient of the torus 
${\bC^{\times}}^n \inj \bC^n$ by the diagonal $\bC^{\times}$.
The action of ${\bC^{\times}}^n$ on $\bC^n$ induces an action of $H$ on $G(k,n)$.
The \emph{Gel'fand--McPherson isomorphism} is a canonical isomorphism
$$U/ \SL(k) \stackrel{\sim}{\map} G^0(k,n)/H.$$
It is obtained by identifying each quotient with the set of generic double $(\SL(k),H)$-orbits on $\bP(M(k,n))$,
where $M(k,n)$ denotes the space of $k \times n$ matrices.
It extends to an isomorphism of Chow quotients
$$({\bP^{k-1}}^{\vee})^n \cq  \SL(k) \stackrel{\sim}{\map} G(k,n)\cq H$$
(\cite{K1}, p. 62, Theorem~2.2.4).
Thus $C_{k,n}$ is canonically identified with $G(k,n)\cq H$.

\subsection{The universal family of cycles over $G(k,n)\cq H$} \label{cycles}

Given a reductive algebraic group $H$ acting on a complex projective variety $X$, one can define
the Hilbert quotient $X \hq H$; the definition is analogous to that of the Chow quotient, with the
appropriate Hilbert scheme playing the r\^{o}le of the Chow variety.
The natural morphism $\Hilb_{\red} \map \Chow$  induces a morphism $X \hq H \map X\cq H$.

In our case, the morphism $G(k,n) \hq H \map G(k,n)\cq H$ is an isomorphism (\cite{K1}, p. 54, Theorem~1.5.2).
The essential point is that the subschemes parametrised by $G(k,n) \hq H$ are reduced. We thus obtain the 
following theorem.

\begin{thm} \label{chowunivfamily}
The universal family of cycles over $G(k,n)\cq H$ is realised by a flat family of reduced subschemes
\begin{eqnarray*}
\begin{array}{ccc}
\cZ & \inj & G(k,n) \times G(k,n)\cq H \\
    & \searrow & \da \\
    &          & G(k,n)\cq H
\end{array}
\end{eqnarray*}
\end{thm}

We pause to describe the Chow quotient of a projective space by a torus in some detail
(Theorem~\ref{chowqtprojspace}, Theorem~\ref{toricunivfamily}).
We use these results to study $G(k,n)\cq H$ via the Plucker embedding.

\begin{notn}
Let $H \cong {\bC^{\times}}^{(n-1)}$ be a torus acting on a projective space $\bP^{N-1}$. 
Pick a linearisation and choose homogeneous
coordinates on $\bP^{N-1}$ so that the action is diagonalised.
Let $A$ be the set of weights of the action in $M:=\Hom(H,\bC^{\times}) \cong \bZ^{n-1}$.
Let $Q$ be the convex hull of $A$ in $M_{\bR}$.
We assume for simplicity that the weights are distinct and that each weight is a vertex of $Q$. 
Let $T={\bC^{\times}}^N / \bC^{\times} \inj \bP^{N-1}$, the big torus.
Then $T$ acts on $\bP^{N-1}$, and the $H$ action is induced by a homomorphism $H \map T$.
We may assume that the action of $H$ on $\bP^{N-1}$ is faithful, then $H \inj T$.
\end{notn}

\begin{defn}
Let $Q$ be a polytope of dimension $d$.
A \emph{polyhedral subdivision} $\underline{Q}$ of $Q$ is a set of polytopes $Q_1,\ldots,Q_r$ of dimension $d$,
such that  
\begin{enumerate}
\item The vertices of $Q_i$ are vertices of $Q$ for each $i$.
\item The intersection $Q_i \cap Q_j$ is a face of $Q_i$ and $Q_j$ for each $i \neq j$.
\item $Q=Q_1 \cup \cdots \cup Q_r$.
\end{enumerate}
A \emph{face} $Q'$ of $\underline{Q}$ is a face of $Q_i$ for some $i$, that is, a cell in
the cellular subdivision of $Q$ defined by $\underline{Q}$.
The subdivision $\underline{Q}$ is \emph{coherent} if there exists a convex piecewise linear function 
$g \colon Q \map \bR$ with domains of linearity given by $\underline{Q}$. 
\end{defn}

The \emph{secondary polytope} $\Sigma Q$ of a polytope $Q$ is a polytope whose faces $F(\underline{Q})$ 
are labelled by coherent polyhedral subdivisions $\underline{Q}$ of $Q$. 
The face $F(\underline{Q})$ is contained in $F(\underline{Q}')$ if and only if $\underline{Q}$ is a 
refinement of $\underline{Q}'$. 
We refer the reader to (\cite{GKZ}, Ch.~7) for more details.

\begin{thm} \cite{KSZ} \label{chowqtprojspace}
The Chow quotient $\bP^{N-1}\cq H$ is a (possibly non-normal) toric variety with torus $T/H$.
Its normalisation is $X_{\Sigma Q}$, the normal toric variety associated to the secondary polytope $\Sigma Q$ of $Q$.
In particular, the toric strata of $\bP^{N-1}\cq H$ are labelled by coherent polyhedral subdivisions of $Q$.
\end{thm}

We review the proof of the theorem.
For $v \in \bP^{N-1}$ generic, we have $\dim \overline{Hv} = n-1$ and $\deg \overline{Hv}=d$, some fixed $d$.
Then 
$$\bP^{N-1}\cq H \inj \Chow_{n-1}(\bP^{N-1},d)$$ 
and 
$$\Chow_{n-1}(\bP^{N-1},d) \inj \bP H^0(\cO_{G(N-n,N)}(d)), \ [Z] \mapsto R_Z,$$
where $R_Z$ is the Chow form of the cycle $Z \subset \bP^{N-1}$.
Now, the $T$-action on $\bP^{N-1}$ induces a $T/H$-action on $\bP^{N-1}\cq H$,
moreover the open embedding $T \inj \bP^{N-1}$ induces an open embedding
$T/H \inj \bP^{N-1}\cq H$; the distinguished point of $\bP^{N-1}\cq H$ is given
by the cycle $Z=\overline{He}$, where $e=(1,\ldots,1) \in \bP^{N-1}$.
Let $B$ be the set of weights of the Chow form $R_Z$
for the $T$-action on $\bP H^0(\cO_{G(N-n,N)}(d))$; they lie in $M_{T/H}=\Hom(T/H,\bC^{\times})$.
Then $\bP\cq H$ is the (possibly non-normal) toric variety $X_B$  
given by $B \subset M_{T/H}$  (as defined in \cite{GKZ}, Ch. 5, p. 166, 1.4).
Let $R$ denote the convex hull of $B$ in $(M_{T/H})_{\bR}$, the \emph{Chow polytope} of $Z$.
Let $X_R$ denote the normal toric variety associated to $R \subset (M_{T/H})_{\bR}$
(\cite{O}, p. 93, Theorem~2.22 or \cite{F}, p. 26--7).
There is a natural equivariant map $X_R \map X_B$ which is finite and birational (it restricts
to the identity on $T/H$), hence $X_R$ is the normalisation of $X_B$.
Finally, the Chow polytope $R$ is the secondary polytope $\Sigma Q$ of $Q$ 
(\cite{KSZ}, p. 209, Theorem~5.1). The proof is rather involved, and proceeds by analysing the extreme toric 
degenerations of $Z=\overline{He} \inj \bP^{N-1}$ corresponding to the vertices of the Chow polytope;
these cycles are sums of coordinate $(n-1)$-planes. Let us at least give a plausible explanation of the result.
Let $X_{\omega}$ denote the homogeneous coordinate on $\bP^{N-1}$ with weight $\omega \in A$. Consider the moment map 
$$\mu \colon \bP^{N-1} \map Q \subset M_{\bR}, \ 
x \mapsto \frac{\sum_{\omega \in A} |X_{\omega}(x)|\cdot \omega}{\sum_{\omega \in A} |X_{\omega}(x)|}.$$ 
A generic orbit closure maps surjectively onto $Q$ under $\mu$.
In general, an orbit closure $\overline{Hv}$ maps onto the weight polytope 
$$\Wt(v)=\conv \{ \omega \in A \ | \ X_{\omega}(v) \neq 0 \}$$
of the point $v$.
A Chow limit $Z$ of generic orbit closures decomposes as a sum $\sum_i Z_i = \sum_i \overline{Hv_i}$
of orbit closures of maximal dimension $n-1$, such that the polytopes $Q_i=\mu(Z_i)$ define a polyhedral decomposition
$\underline{Q}$ of $Q$. The toric stratum of $\bP^{N-1}\cq H$ containing $[Z]$ is thus labelled by $\underline{Q}$.

The description of Chow limits of orbit closures via the moment map is unfortunately not precise enough
for our purposes. We obtain a more careful description below by constructing a universal family
over $X_{\Sigma Q}$, the normalisation of $\bP^{N-1}\cq H$. We first need a definition.
 
\begin{defn} (cf. \cite{A3}, p. 28, Defn.~2.4.2) \label{degtoricdef}
Suppose given a polyhedral subdivision $\underline{Q}$ of a lattice polytope $Q \subset M_{\bR}$.
We define the \emph{degenerate toric variety} $X_{\underline{Q}}$ associated to $\underline{Q}$ as follows.
Let $M'=\bZ \oplus M$ and embed $Q$ in $M'_{\bR}$ by identifying $M_{\bR}$ with the affine hyperplane $(1,M_{\bR})$.
Let $\Omega$ be the fan of cones over the faces of $\underline{Q}$ in $M'_{\bR}$.
Write $S=|\Omega| \cap M'= \Cone(Q) \cap M'$.
Let $R[\Omega]$ denote the $\bC$-algebra with $\bC$-basis $\{ \chi^s \ | \ s \in S \}$ and
multiplication
\begin{equation*}
\chi^s \cdot \chi^t = \left\{ \begin{array}{ll}
                              \chi^{s+t} & \mbox{if there is a cone of $\Omega$ containing $s$ and $t$} \\
                              0          & \mbox{otherwise}
\end{array} \right.
\end{equation*}
We define $X_{\underline{Q}}=\Proj R[\Omega]$.
\end{defn}

We note immediately that $X_{\underline{Q}}$ has normal toric components 
$$X_{Q_i}=\Proj(\bC[\Cone(Q_i) \cap M'])$$ 
corresponding to the maximal polytopes $Q_i$ of the subdivision $\underline{Q}$, and these are glued together following the 
same combinatorial rules used to glue the $Q_i$ to form $Q$. Moreover, the sheaf $\cO_{X_{\underline Q}}(1)$
is a line bundle obtained by gluing the line bundles $\cO_{X_{Q_i}}(1)$ on the components,
and $H^0(\cO_{X_{\underline{Q}}}(1))$ has a canonical basis identified with the set of lattice points of $Q$.

\begin{thm} (cf. \cite{L}, p.~15--16) \label{toricunivfamily}
There is a normal toric variety U, with torus $T$, and an equivariant finite map 
$U \map \bP^{N-1} \times X_{\Sigma Q}$, such that the image of $U$ in $\bP^{N-1} \times X_{\Sigma Q}$
is the pullback of the universal family of cycles over $\bP^{N-1}\cq H$. The reduced fibres of $U / X_{\Sigma Q}$
over the orbit labelled by a polyhedral subdivision $\underline{Q}$ of $Q$
are isomorphic to the degenerate toric variety $X_{\underline{Q}}$ associated to $\underline{Q}$.
\end{thm}

\begin{proof}
The family $U \map X_{\Sigma Q}$ is defined explicitly by a map of fans (cf. \cite{L}, p.137--145)
as follows. Let $\bR^A$ denote the set of functions $\psi \colon A \map \bR$.
For $\psi \in \bR^A$, let $g_{\psi} \colon Q \map \bR$ denote the 
convex piecewise-linear function whose graph is the lower envelope of the polytope
$$\conv \{(\omega, \psi(\omega)) \ | \ \omega \in A \} \subset M_{\bR} \times \bR.$$
For $\underline{Q}$ a subdivision of $Q$, define 
$$C(\underline{Q})=\{ \psi \in \bR^A \ | \ g_{\psi} \mbox{ is linear over each polytope in $\underline{Q}$} \},$$
a closed cone in $\bR^A$.
Let $\Aff(A) \inj \bR^A$ denote the subspace of affine linear functions $\psi$, i.e., functions $\psi$ such that 
$g_{\psi}$ is (affine) linear over the whole of $Q$.
Thus $\Aff(A) \subset C(\underline{Q})$ for each $\underline{Q}$.
The quotient cones
$$C(\underline{Q}) / \Aff(A) \subset \bR^A / \Aff(A)$$
constitute a complete fan (cf. \cite{L}, p.~52, Proposition~II.1), which we denote by $\Delta(Q)$.
The vector space $\bR^A / \Aff(A)$ contains the lattice $\bZ^A / \Aff(A) \cap \bZ^A$ given by $\bZ$-valued functions,
which is naturally identified with the dual $N_{T/H}$ of the characters $M_{T/H}=\Hom(T/H,\bC^{\times})$ of $T/H$.
The fan $\Delta(Q)$ together with this lattice defines a complete toric variety with torus $T/H$.
In fact, this toric variety is the projective toric variety $X_{\Sigma Q}$ defined by the secondary polytope $\Sigma Q$,
since $\Delta(Q)$ coincides with the normal fan of $\Sigma Q \subset (M_{T/H})_{\bR}$ (\cite{GKZ}, Theorem~2.4, p.~228). 

For $\underline{Q}$ a subdivision of $Q$ and $Q'$ a face of $\underline{Q}$, define 
$$C(\underline{Q},Q')=\{ \psi \in C(\underline{Q}) \ | \ g_{\psi} \mbox{ attains its minimum on $Q'$} \},$$
a closed cone in $\bR^A$. Let $\bR \inj \bR^A$ denote the constant functions.
The quotient cones
$$C(\underline{Q},Q')/ \bR \inj \bR^A / \bR$$
form a complete fan (cf. \cite{L}, p.~137, Proposition~IV.1) which we denote $\tilde{\Delta}(Q)$.
The lattice $\bZ^A / \bZ \inj \bR^A / \bR$ is naturally identified with the cocharacters $N_T$ of $T$.
The fan $\tilde{\Delta}(Q)$ together with this lattice defines a complete toric variety $U$ with torus $T$.
The equivariant morphism $U \map X_{\Sigma Q}$ is defined by 
the map $(N_T,\tilde{\Delta}(Q)) \map (N_{T/H},\Delta(Q))$ of fans induced by the quotient map 
$\bR^A/ \bR \map \bR^A / \Aff(A)$.

The reduced fibres of the map $f \colon U \map X_{\Sigma Q}$ may be described by analysing the map of fans.
Using the result and notation of Lemma~\ref{toricfibres},  
$$f^{-1}O_{C(\underline{Q})}=\coprod_{Q'} O_{C(\underline{Q},Q')}.$$ 
So, the fibres of $f$ over the orbit on $X_{\Sigma Q}$ labelled by $\underline{Q}$ have strata labelled by
the faces $Q'$ of $\underline{Q}$. One checks that the strata are precisely the toric varieties $X_{Q'}$ defined by
$Q'$, and the gluing is as described in Definition~\ref{degtoricdef}. Hence the fibres are isomorphic to the degenerate
toric variety $X_{\Sigma Q}$, as desired.

The fan $\tilde{\Delta}(Q) \inj \bR^A / \bR$ is a refinement of the fan corresponding to $\bP^{N-1}$,
that is, the fan generated by the standard basis of $\bR^A$. We thus obtain an equivariant birational map
$U \map \bP^{N-1}$. On the fibre $X_{\underline{Q}}$ of $U/X_{\Sigma Q}$ 
over the distinguished point of the stratum labelled by $\underline{Q}$, the map 
$X_{\underline{Q}} \map \bP^{N-1}$ is defined by the global sections of $\cO_{X_{\underline{Q}}}(1)$ corresponding 
to the vertices $A$ of $Q$. The map on an arbitrary fibre is defined similiarly, by $T$-equivariance.
In particular, the map $g \colon U  \map \bP^{N-1} \times X_{\Sigma Q}$ is finite.
Let $V$ denote the universal family over the Chow quotient $\bP^{N-1}\cq H$; it remains to show that
$g_{\star}U$ is the pullback of $V$. 
The image of the fibre of $U$ over the distinguished point of $X_{\Sigma Q}$ is $\overline{He} \inj \bP^{N-1}$,
the fibre of $V$ over the distinguished point of $\bP^{N-1}\cq H$.
Hence $g_{\star}U$ agrees with the pullback of $V$ over $T/H$ by $T$-equivariance; it follows that $g_{\star}U$ 
is equal to the pullback of $V$.
\end{proof}

\begin{lem} \label{toricfibres}
Let $f \colon X' \map X$ be a equivariant map of toric varieties given by
a map $\phi \colon (N',\Delta') \map (N,\Delta)$ of fans.
Let $O_{\tau}$ denote the orbit corresponding to a cone $\tau$.  
Then, for $\sigma \in \Delta$,
$$f^{-1}(O_{\sigma})= \coprod_{\sigma' \in S} O_{\sigma'},$$
where 
$S=\{ \sigma' \in \Delta' \ | \mbox{ $\sigma$ is the smallest cone of $\Delta$ containing $\phi(\sigma')$} \}.$
\end{lem}

We can now describe the fibres of the universal family $\cZ / (G(k,n)\cq H)$ carefully.
The Plucker embedding $G(k,n) \inj \bP$ is $H$-equivariant. 
It thus yields an embedding of Chow quotients $G(k,n)\cq H \inj \bP\cq H$, such that the 
universal family of cycles over $G(k,n)\cq H$ is the pullback of the universal family over $\bP\cq H$. 
The weight polytope $Q$ of the $H$-action on $\bP=\bP(\wedge^k \bC^n)$ is
the $(k,n)$-\emph{hypersimplex}   
$$\Delta(k,n)=\conv\{ e_I \ | \ |I|=k, \ I \subset [n] \} \subset ( \sum_1^{n} x_i = k ) \subset \bR^n.$$
Here $e_I= \sum_{i \in I} e_i$, the weight of the Plucker coordinate $P_I$. 
The affine space $(\sum_1^{n} x_i = k)$ is identified with
$M_{\bR}= \Hom(H,\bC^{\times})_{\bR} = (\sum_1^{n} x_i = 0)$ by fixing an origin.

\begin{thm} \label{chowqtgrass}
The Chow quotient $G(k,n)\cq H$ has a stratification with strata labelled by polyhedral subdivisions 
$\underline{Q}$ of $\Delta(k,n)$ such that each edge of $\underline{Q}$ is also an edge of $\Delta(k,n)$.
Let $\cZ/(G(k,n)\cq H)$ denote the universal family of cycles.
The fibres of $\cZ$ over the locally closed stratum corresponding to a subdivision $\underline{Q}$ are 
isomorphic to the degenerate toric variety $X_{\underline{Q}}$ associated to $\underline{Q}$. 
More precisely, let $(G(k,n)\cq H)'$ denote the strict transform of $G(k,n)\cq H$ in the normalisation
$X_{\Sigma Q}$ of $\bP\cq H$, and $\cZ'$ the pullback of $\cZ$.
Then $\cZ'/(G(k,n)\cq H)'$ is identified with the restriction of the family $U/X_{\Sigma Q}$
of \ref{toricunivfamily}.
\end{thm} 

\begin{proof}
The weight polytope of a point $[L] \in G(k,n)$ is a matroid polytope (\cite{GS}, p.~143--145).
Concretely, a polytope $Q' \subset \Delta(k,n)$ is a matroid polytope if each edge of $Q'$ is also an edge of
$\Delta(k,n)$ (op. cit., p.~144, Theorem~1). 
We say a polyhedral subdivision $\underline{Q}$ of $\Delta(k,n)$ is a matroid decomposition if each polytope 
in $\underline{Q}$ is a matroid polytope. 
Let $(\bP\cq H)^0$ denote the union of the toric strata of $\bP\cq H$ corresponding to 
matroid decompositions. Then $(\bP\cq H)^0$ is an open subvariety of $\bP\cq H$ and contains $G(k,n)\cq H$ by construction. 
Let $(X_{\Sigma Q})^0$ denote the corresponding open subvariety of the normalisation $X_{\Sigma Q}$ of $\bP\cq H$
and $U^0/ (X_{\Sigma Q})^0$ the pullback of the family $U/ X_{\Sigma Q}$. 
Then $U^0/(X_{\Sigma Q})^0$ has reduced fibres and the natural map $U^0 \map \bP \times (X_{\Sigma Q})^0$ is a 
closed embedding  (\cite{L}, p. 141, Proposition~IV.3).
It follows that  $\cZ' / (G(k,n)\cq H)'$ is identified with the restriction of
$U/X_{\Sigma Q}$ by Theorem~\ref{toricunivfamily}.
\end{proof}

\subsection{The universal family of stable pairs} \label{slice}

We obtain the universal family of stable pairs over $C_{k,n}=G(k,n)\cq H$ as a transversal slice of the 
universal family of cycles $\cZ / (G(k,n)\cq H)$.
The same construction has already been applied to a generic fibre of $\cZ/C_{k,n}$ by Kapranov;
in this case one recovers the corresponding hyperplane arrangement.

We first describe the construction on fibres.
A fibre $Z$ of $\cZ/C_{k,n}$ is a degenerate toric variety corresponding to a polyhedral subdivision 
of $\Delta(k,n)$. Let $B$ denote the toric boundary of $Z$, i.e., the divisor on $Z$ corresponding to the boundary
of $\Delta(k,n)$. 
The facets of $\Delta(k,n)$ are 
$$\Gamma_i^+=\conv\{e_I \ | \ i \in I \}, \ i=1,\ldots,n, $$
$$\Gamma_i^-=\conv\{e_I \ | \ i \notin I \}, \ i=1,\ldots,n. $$ 
Let $B=\sum B_i^+ + \sum B_i^-$  be the corresponding decomposition of $B$.
Then 
$$B_i^+=Z \cap G(k-1,n-1)_{e_i},$$
$$B_i^-=Z \cap G(k,n-1)_i$$
(as schemes), where
$$G(k-1,n-1)_{e_i}=\{L \in G(k,n) \ | \ e_i \in L \} \inj G(k,n),$$
$$G(k,n-1)_i=\{L \in G(k,n) \ | \ L \subset (x_i=0) \} \inj G(k,n)$$
(\cite{K1}, p. 56, Proposition~1.6.10).
Let $e= (1,\ldots,1) \in \bC^n$ and write
$$G(k-1,n-1)_e =\{L \in G(k,n) \ | \ e \in L \} \inj G(k,n).$$

\begin{defn}
Given $[Z] \in G(k,n)\cq H$, we define the associated pair
$(X,D=D_1+\cdots+D_n)$ by 
\begin{eqnarray*}
\begin{array}{ccc}
X & = & Z \cap G(k-1,n-1)_e, \\ 
D_i & = & B^+_i \cap G(k-1,n-1)_e.
\end{array}
\end{eqnarray*}
\end{defn}
\begin{rem} Note immediately that $B_i^- \cap G(k-1,n-1)_e = \emptyset$,
hence $D=B|_X$.
\end{rem}

We eventually prove that the pair $(X,D)$ associated to $Z$ is a stable pair.
We first establish that $(X,D)$ is an slc pair of dimension $k-1$. 
We require the following transversality result.

\begin{lem} \label{transversality}
Consider the map 
$$f \colon H \times G(k-1,n-1)_e \map G(k,n), \ (h,[L]) \mapsto h[L]$$
obtained from the $H$-action on $G(k,n)$.
The map $f$ is smooth of relative dimension $k-1$, with  image $G(k,n)' = G(k,n) \backslash  \bigcup G(k,n-1)_i$, the locus
of $k$-planes not contained in a coordinate hyperplane.
\end{lem}
\begin{proof}
Let
\begin{eqnarray*}
\begin{array}{ccc}
S & \inj & \bP^{n-1} \times G(k,n) \\
  &      &  \da p\\
  &      &  G(k,n)
\end{array}
\end{eqnarray*}
denote the projectivised tautological bundle over $G(k,n)$.
Let $H = {\bC^{\times}}^n/\bC^{\times} \inj \bP^{n-1}$ be the usual torus embedding.
Then $f$ factorises as $p \circ g$, where
$$g \colon H \times G(k-1,n-1)_e  \stackrel{\sim}{\map} (H \times G(k,n)) \cap S \inj \bP^{n-1} \times G(k,n), $$
$$(h,[L]) \mapsto (h, h[L]).$$ Hence $f$ is smooth, 
and the fibre of $f$ over $[L] \in G(k,n)$ is identified with $\bP(L) \cap H$.
\end{proof}

\begin{prop} \label{(X,D)slc}
The pair $(X,D)$ is slc of dimension $k-1$.
More precisely, the map $f_{Z} \colon H \times X \map Z$ induced by the $H$-action is smooth.
Thus, everywhere locally on $X$, the pair $(X,D)$ is obtained as a transversal slice of $(Z,B)$ 
at a $H$-orbit of codimension $\le k-1$.
\end{prop}
\begin{proof}
The map $f_{Z}$ is obtained from the map $f$ of 
Lemma~\ref{transversality} by the base change $Z \inj G(k,n)$; thus $f_{Z}$
is smooth of relative dimension $k-1$, with image $Z \cap G(k,n)'= Z \backslash \bigcup B_i^-$.
The degenerate toric pair $(Z,B)$ is slc by a well-known toric argument (cf. \cite{A2}, p.~9, Lemma~3.1).
Hence $(X,D)$ is slc since $f_Z$ is smooth and $f_Z^{\star}B=H \times D$.
\end{proof}

The scheme $X$ has a stratification with strata given by arbitrary intersections of the 
irreducible components of $X$ and the divisor $D$. The associated poset captures the
combinatorial type of the pair $(X,D)$. We use the embedding $X \inj Z$ to
describe this stratification.

\begin{thm} \label{Xstrata}
Let $(X,D)$ be the pair associated to $[Z] \in G(k,n)\cq H$. Consider the stratification of $X$ defined by
the irreducible components of $X$ and $D$.
\begin{enumerate}
\item The stratification coincides with that induced by the stratification of $Z$ by orbit closures.
The non-empty strata of $X$ are the intersections $Y_e= Y \cap G(k-1,n-1)_e$, where $Y$ is a 
stratum of $Z$ not contained in $\cup B^-_i$.
Each stratum $Y_e$ is an irreducible, normal variety of the expected dimension $\dim Y - (n-k)$.
Moreover, $Y_e$ is rational and has rational singularities. 
\item Let $\underline{Q}$ denote the polyhedral decomposition of $Q=\Delta(k,n)$ corresponding to the degenerate toric variety $Z$.
Then the poset of strata of $X$ is identified with the poset of faces of $\underline{Q}$ which are not contained in $\cup \Gamma_i^-$.
\end{enumerate}
\end{thm}
\begin{proof}
Let $Y$ be a toric stratum of $Z$.
Assume that $Y$ is not contained in $\bigcup B^-_i$ (otherwise $Y \cap G(k-1,n-1)_e = \emptyset$).
Again consider the map $f_{Y} \colon H \times Y_e \map Y$ induced by the $H$-action;
it is smooth of relative dimension $k-1$, with image $Y \cap G(k,n)'$.
Thus $Y_e$ is irreducible of dimension $\dim Y - (n-k)$.
It follows that the stratification of $X$ is induced by the toric stratification of $Z$ as claimed.
Moreover, $Y_e$ is normal and has rational singularities since $Y$ satisfies these 
properties (the stratum $Y$ is a normal toric variety).
Finally, $Y_e$ is rational by \cite{K1}, p.~68, Proposition~3.1.9(a).
\end{proof}

\begin{rem}
The combinatorial type of $(X,D)$ may be alternatively described by the \emph{dual complex} of $(X,D)$.
The dual complex of $(X,D)$ is a pair $(\Sigma, \partial \Sigma)$ of CW-complexes.
For each stratum $Y$ of $X$ there is a cell $\sigma_Y$ of $\Sigma$ of dimension $\dim X - \dim Y$,
such that $\sigma_Y \subset \sigma_{Y'}$ if $Y' \subset Y$.
For each stratum $Y$ contained in $D$, there is in addition a cell $\sigma^{\partial}_Y$ of 
$\partial \Sigma$ of dimension $\dim D - \dim Y$, such that $\sigma^{\partial}_Y \subset \sigma_Y$. 
For example, if $k=2$, the dual complex is the dual graph of the stable curve $(X,D)$,
a tree with $n$ labelled endpoints.
The complex $\Sigma$ is contractible, and
$\partial \Sigma$ is homeomorphic to the $(k-2)$-skeleton of an $(n-1)$-simplex.
The subdivision $\underline{Q}$ determines the dual complex by \ref{Xstrata}.
Conversely, one shows that the homeomorphism type of the dual complex determines the
subdivision $\underline{Q}$. We omit the details.
\end{rem}

The pairs $(X,D)$ fit into a flat family $(\cX^u,\sD^u=\sD^u_1+ \cdots+ \sD^u_n)/C_{k,n}$ defined by
repeating the construction in the relative context. That is,
\begin{eqnarray*}
\begin{array}{ccc}
\cX^u & = & \cZ \cap (G(k-1,n-1)_e \times C_{k,n}) \\ 
\sD_i & = & \cB^+_i \cap (G(k-1,n-1)_e \times C_{k,n})
\end{array}
\end{eqnarray*}
where
$$\cB_i^+=\cZ \cap (G(k-1,n-1)_{e_i} \times C_{k,n}).$$

As remarked above, the generic fibres $(X,D)$ were identified by Kapranov.
We record his precise result here for future reference.
\begin{prop}(\cite{K1}, p. 71, Proposition~3.2.3 and p. 75, Proposition~3.4.3) \label{genericfibres}
Let $M_1,\ldots,M_n \subset \bP^{k-1}$ be an arrangement of hyperplanes in general position.
Let $[Z] \in G(k,n)\cq H$ be the point corresponding to this arrangement 
under the Gel'fand--McPherson isomorphism, and $(X,D)$ the associated pair.
Then $(X,D)$ is canonically isomorphic to $(\bP^{k-1},M_1+\cdots+M_n)$.
Moreover, the embedding $(X,D) \inj G(k-1,n-1)_e$ is given by the vector bundle 
$\Omega_{\bP^{k-1}}(\log M)$ of differential forms on $\bP^{k-1}$
with logarithmic poles along the hyperplanes $M_i$, together with the 
isomorphism
$$H^0(\Omega_{\bP^{k-1}}(\log M)) \stackrel{\sim}{\map} \ker (\bC^n \stackrel{\Sigma}{\map} \bC)$$
given by taking residues along the $M_i$.
\end{prop}

\subsection{Restriction to a hyperplane}
Given a generic hyperplane arrangement, the intersection of the $i$th hyperplane with the remaining hyperplanes is 
again a hyperplane arrangement. We thus obtain maps $a_i \colon H_{k,n}  \map  H_{k-1,n-1}$ given by 
$$(\bP^{k-1},M_1+\cdots+M_n)  \mapsto  (M_i, (M_1+\cdots+M_{i-1}+M_{i+1}+ \cdots+M_n)|_{M_i}).$$
These maps extend to maps $\tilde{a}_i \colon C_{k,n} \map C_{k-1,n-1}$ of the Chow quotient compactifications,
given by
$$G(k,n)\cq H  \map  G(k-1,n-1)_{e_i}\cq H_i, \ \  [Z]  \mapsto  [Z \cap G(k-1,n-1)_{e_i}]$$ 
(\cite{K1}, p.~56, Theorem~1.6.6).  
Here $H_i$ is the quotient of $H={\bC^{\times}}^n/\bC^{\times}$ by the $i$th coordinate $\bC^{\times}$.
We observe that, using stable pairs, 
the extension $\tilde{a}_i$ of $a_i$ may be described in the same way as $a_i$.

\begin{prop} \label{restriction}
Let $(X,D)$ be the pair associated to $[Z] \in G(k,n)\cq H$.
Then $(D_i,(D-D_i)|_{D_i})$ is the pair associated to $\tilde{a}_i([Z])$.
\end{prop}
\begin{proof}
Recall that the pair $(X,D)$ is defined by 
$X=Z\cap G(k-1,n-1)_e$ and $D_i=B^+_i \cap G(k-1,n-1)_e$ where 
$B_i^+ = Z \cap G(k-1,n-1)_{e_i}$.
The result follows easily.
\end{proof}

\section{The sheaf of logarithmic differentials on a stable pair}

Let $(X,D)$ be the pair associated to some $[Z] \in G(k,n)\cq H$.
We show that the embedding $(X,D) \inj G(k-1,n-1)_e$ 
may be defined using a certain sheaf of logarithmic differentials on 
$(X,D)$ (cf. Proposition~\ref{genericfibres}), and use this to prove
that the embedding is canonical. In order to define this sheaf, we must 
first define a log structure on $(X,D)$ (since $X$ may be reducible).

\subsection{Definition of the log structure on $(X,D)$}\label{logstr}

We assume that the reader is familiar with log structures (see, e.g., \cite{KK}).
Recall that $X= Z \cap G(k-1,n-1)_e$. The scheme $Z$ is a fibre of the universal family 
$U \map X_{\Sigma Q}$ of Theorem~\ref{toricunivfamily}, where $Q=\Delta(k,n)$.
The varieties $U$ and $X_{\Sigma Q}$ are toric with tori $T$ and $T/H$, and the map
$U \map X_{\Sigma Q}$ is equivariant with respect to $T \map T/H$. Hence $U$ and $X_{\Sigma Q}$ have canonical log 
structures $U^{\dagger}$ and $X_{\Sigma Q}^{\dagger}$, 
and the map $U \map X_{\Sigma Q}$ extends naturally to a log smooth map $U^{\dagger} \map X_{\Sigma Q}^{\dagger}$
of log schemes. We define the log structure $Z^{\dagger}/k^{\dagger}$ on the fibre $Z/k$ by base change.
Finally, we define the log structure $X^{\dagger}/k^{\dagger}$ on $X /k$ by restriction.

\subsection{Identification of the sheaf of log differentials} \label{logdiff}

We describe the sheaf of log differentials on $X/k$ defined by the log structure.
Roughly speaking, it is the sheaf of differentials on $X$ with logarithmic poles along
the divisor $D$ and the double locus, such that the two residues along a component 
of the double locus sum to zero. The log structure is required to glue the component 
sheaves at the double locus.

For $Y$ a toric variety with torus $T$, let 
$\Omega_{Y^{\dagger}/k}$ denote the sheaf of differentials with logarithmic poles along the toric boundary.
Equivalently, $\Omega_{Y^{\dagger}/k}$ is the sheaf of log differentials associated to the 
canonical log structure on $Y$. There is a canonical isomorphism
$$\cO_Y \otimes_{\bZ} M_T \stackrel{\sim}{\map} \Omega_{Y^{\dagger}/k}, \ m \mapsto d\chi^m/\chi^m$$
where $M_T=\Hom(T,\bC^{\times})$, and, for $m \in M_T$, $\chi^m$ denotes the corresponding
rational function on $Y$ (\cite{F}, p.~87). We establish an analogous result for the degenerate 
toric variety $Z$.

\begin{lem}
There is a canonical isomorphism
$$\cO_Z \otimes_{\bZ} M \stackrel{\sim}{\map} \Omega_{Z^{\dagger}/k^{\dagger}}$$ 
The sheaf $\Omega_{Z^{\dagger}/k^{\dagger}}$ can be obtained by gluing the sheaves of log differentials on
the components $Z_i$ of $Z$ using the isomorphisms $\Omega_{Z_i^{\dagger}/k} \cong \cO_{Z_i} \otimes_{\bZ} M.$
\end{lem}

\begin{proof}
Recall that $Z$ is a fibre of the equivariant morphism of toric varieties $U \map X_{\Sigma Q}$.
The exact sequence
$$0 \map f^{\star} \Omega_{X_{\Sigma Q}^{\dagger}/k} \map \Omega_{U^{\dagger}/k} \map 
\Omega_{U^{\dagger}/X_{\Sigma Q}^{\dagger}} \map 0$$
yields
$\Omega_{U^{\dagger}/ X_{\Sigma Q}^{\dagger}} \cong \cO_U \otimes_{\bZ} M,$
and thus 
$\Omega_{Z^{\dagger}/k^{\dagger}} \cong \cO_Z \otimes_{\bZ} M$
by base change.
\end{proof}

\begin{notn}
Let $d\chi^m / \chi^m$ denote the global section of $\Omega_{Z^{\dagger}/k^{\dagger}}$
corresponding to $m \in M$. It restricts to the usual log differential $d\chi^m / \chi^m$
on each component $Z_i$ of $Z$. Note that $\chi^m$ is \emph{not} a rational function on $Z$,
i.e., it is not an element of the total quotient ring of $\cO_Z$, although it is of course 
a well defined rational function on each component $Z_i$. 
\end{notn}
Write $G=G(k,n)$, and let $\cS$ and $\cQ$ denote the tautological subsheaf and quotient sheaf on $G$.
Also write $G_e=G(k-1,n-1)_e$, and define $\cS_e$ and $\cQ_e$ similiarly.
Let $\mathfrak{h}$ denote the Lie algebra of $H={\bC^{\times}}^n/\bC^{\times}$. 
Note that the ambient space $\bC^n/\bC e$ for the Grassmannian $G_e$ is identified with $\mathfrak{h}$.
We write $M=\Hom(H,\bC^{\times})$ as before, then $M_{\bC}=\mathfrak{h}^{\vee}$.

\begin{thm} \label{logdiffX}
\begin{enumerate}
\item The embedding $X \inj Z$ induces an isomorphism $$\cS^{\vee}_e|_X \stackrel{\sim}{\map} \Omega_{X^{\dagger}/k^{\dagger}}.$$ 
The corresponding map $\mathfrak{h}^{\vee} \map H^0(\Omega_{X^{\dagger}/k^{\dagger}})$ is given by 
$m \mapsto (d\chi^m / \chi^m )|_X$ for $m \in M \subset \mathfrak{h}^{\vee}$.
\item The sheaf $\omega_X(D)$ is identified with the top exterior power $\wedge^{k-1} \Omega_{X^{\dagger}/k^{\dagger}}$
of the sheaf of logarithmic differentials. Thus $\omega_X(D) \cong \cO_{G_e}(1)|_X$ 
where $\cO_{G_e}(1)$ denotes the Plucker line bundle on $G_e$. In particular, $\omega_X(D)$ is a very ample line bundle.\end{enumerate}
\end{thm}

\begin{proof}
We compute $\Omega_{X^{\dagger}/k^{\dagger}}$ using the exact sequence
\begin{equation} \label{OmegaX}
0 \map \cI/\cI^2 \map \Omega_{Z^{\dagger}/k^{\dagger}}|_X \map \Omega_{X^{\dagger}/k^{\dagger}} \map 0,
\end{equation}
where $\cI=\cI_{X/Z}$ is the ideal sheaf of $X$ in $Z$. 
Recall that $X = Z \cap G_e$.
The embedding $G_e \inj G$ is a local complete intersection defined by the vanishing of the section of $\cQ$
given by $e \in \bC^n$. In particular, the normal bundle $\cN_{G_e/G}$ is identified with $\cQ|_{G_e}=\cQ_e$.
Restricting to $Z$, we deduce that the embedding $X \inj Z$ is also a local complete intersection, with normal
bundle $\cN_{X/Z}=\cQ_e|_X$ (note that $X$ has the expected dimension by Proposition~\ref{(X,D)slc}).
The sheaf 
$$T_{Z^{\dagger}/k^{\dagger}}=\Omega_{Z^{\dagger}/k^{\dagger}}^{\vee} \cong \cO_Z \otimes_{\bC} \mathfrak{h}$$
of vector fields on $Z$ with logarithmic zeroes is the sheaf generated by the vector fields induced
by the $H$-action.
The map $T_{Z^{\dagger}/k^{\dagger}} \map \cN_{X/Z}$ gives the corresponding first order deformations of $X \inj Z$.
One checks that this map is identified with the restriction of the quotient map 
$\cO_{G_e} \otimes \mathfrak{h} \map \cQ_e$ on $G_e=G(k-1, \mathfrak{h})$.
Hence the exact sequence (\ref{OmegaX}) is identified with the dual of the exact sequence
$$0 \map \cS_e \map \cO_{G_e} \otimes_{\bC} \mathfrak{h} \map \cQ_e \map 0$$
restricted to $X$.
In particular, the sheaf $\Omega_{X^{\dagger}/k^{\dagger}}$ is identified with $\cS_e^{\vee}|_X$.

The top exterior power $\wedge^{k-1}\Omega_{X^{\dagger}/k^{\dagger}}$ is thus identified with the 
restriction of the Plucker line bundle $\cO_{G_e}(1)$. On the other hand, we have 
$\wedge^{n-1}\Omega_{Z^{\dagger}/k^{\dagger}} = \omega_Z(B)$, where $B$ denotes the toric boundary of $Z$.
We deduce from the exact sequence (\ref{OmegaX}) that $\wedge^{k-1}\Omega_{X^{\dagger}/k^{\dagger}}=\omega_X(D)$, 
using adjunction for the dualising sheaf and the equality $D=B|_X$.
\end{proof}

\subsection{The canonical basis of $H^0(\omega_X(D))$} \label{sectcanbasis}

Our aim in this section is to prove that the embedding $(X,D) \inj G_e$ is canonically determined by $(X,D)$.
The embedding is given by the locally free sheaf $\Omega_{X^{\dagger}/k^{\dagger}}$ defined
above, together with a certain map $\mathfrak{h}^{\vee} \map H^0(\Omega_{X^{\dagger}/k^{\dagger}})$.
Equivalently, after composing with the Plucker embedding, the embedding is given by the line bundle
$\omega_X(D) = \wedge^{k-1}~\Omega_{X^{\dagger}/k^{\dagger}}$, together with the induced map
$\wedge^{k-1}~\mathfrak{h}^{\vee} \map H^0(\omega_X(D))$. We prove that the embedding is canonical 
by identifying this map.  

\begin{rem}
We expect that the map $\mathfrak{h}^{\vee} \map H^0(\Omega_{X^{\dagger}/k^{\dagger}})$ is the inverse of an
isomorphism 
$$H^0(\Omega_{X^{\dagger}/k^{\dagger}}) \map \mathfrak{h}^{\vee} = \ker ( \bC^n \stackrel{\Sigma}{\map} \bC)$$
defined by taking residues along the divisors $D_1,\ldots,D_n$ (cf. Proposition~\ref{genericfibres}). However,
I don't know how to prove this at present. It is also unclear a priori that the sheaf $\Omega_{X^{\dagger}/k^{\dagger}}$
is canonically determined by $(X,D)$. So, in the following, we work instead with the sheaf $\omega_X(D)$.
\end{rem}

\begin{lem} \label{points}
The intersection $D_{i_1} \cap \cdots \cap D_{i_{k-1}}$ is a point for each $i_1< \cdots <i_{k-1}$.
At each such point $P=P_{i_1,\ldots,i_{k-1}}$, the scheme $X$ is smooth and the divisor 
$D$ has normal crossings. Thus
$$(P \in X,D) \cong (0 \in \bC^{k-1}, (x_1\cdots x_{k-1}=0))$$ 
where $D_{i_j}$ is identified with $(x_j=0)$.
\end{lem}

\begin{proof}
Write $J=\{ i_1,\cdots,i_{k-1} \}$.
The subscheme $D_J = \bigcap_{i \in J} D_i \subset X$ is the intersection of the subscheme 
$B^+_J =\bigcap_{i \in J} B^+_i \subset Z$ with $X$. The scheme $B_J$ corresponds to the face 
$\Gamma= \bigcap_{i \in J} \Gamma^+_i$ of $\Delta(k,n)$. The face $\Gamma$ is the $(n-k)$-simplex
$$\conv \{ e_{J\cup \{ i \}} \ | \ i \notin J \}.$$
In particular, $\Gamma$ is a face of the subdivision $\underline{Q}$ of $Q=\Delta(k,n)$ corresponding to $Z$
(since the only subdivision of a simplex is the trivial one).
Thus $B_{J}^+$ is a single toric stratum of $Z$, of dimension $n-k$. By Theorem~\ref{Xstrata}, the corresponding
stratum $D_J$ of $X$ is a (reduced) point.

To prove the second part, we analyse the subdivision $\underline{Q}$ at $\Gamma$. Let $e_I$ be a vertex of $\Gamma$,
so $I=J \cup \{ i_k \}$, some $i_k$, and fix the embedding $Q \inj M_{\bR}$ by 
identifying $e_I$ as the origin.
Let $\langle S \rangle$ denote the cone and $\langle S \rangle_{\bR}$ the vector space generated by a set 
$S \subset M_{\bR}$. Consider the quotient cone
$$\overline{\Cone}_{\Gamma}(Q) = ( \langle Q \rangle + \langle \Gamma \rangle_{\bR} ) / 
\langle \Gamma \rangle_{\bR}.$$
We have 
$$\langle Q \rangle = \langle e_j-e_i \ | \ j \notin I, \ i \in I \rangle$$
$$\langle \Gamma \rangle = \langle e_j -e_{i_k} \ | \ j \notin I \rangle.$$
So, identifying $M_{\bR}/ \langle \Gamma \rangle_{\bR}$ with
$(x_j=0, j \notin I) \subset M_{\bR}$, 
$$\overline{\Cone}_{\Gamma}(Q)=\langle e_{i_k} - e_i \ | \ i \in J \rangle.$$ 
In particular, $\overline{\Cone}_{\Gamma}(Q)$ is simplicial, and the generators of
$\langle Q \rangle$ yield a minimal set of generators of $\overline{\Cone}_{\Gamma}(Q)$.
It follows that there is a unique maximal polytope $Q'$ of $\underline{Q}$ containing
$\Gamma$. For, the edges  of any such $Q'$ are also edges of $Q$ (Theorem~\ref{chowqtgrass}),
whence $\overline{\Cone_{\Gamma}}(Q')=\overline{\Cone_{\Gamma}}(Q)$. 
Hence $Z$ has a unique component $Z'$
containing the stratum $B^+_J$ corresponding to $\Gamma$, and $Z'$ is smooth at the generic point of $B^+_J$
(since $\overline{\Cone_{\Gamma}}(Q)$  is simplicial). By smoothness of $f_Z \colon H \times X \map Z$ 
(Proposition~\ref{(X,D)slc}) we deduce that $X$ is smooth at $D_J$.
\end{proof}

The Poincar\'{e} residue map
$$\Omega_{Z^{\dagger}/k^{\dagger}} \map (\oplus_i \, \cO_{B_i^+}) \oplus (\oplus_i \, \cO_{B_i^-})$$
(cf. \cite{F}, p.~87, \cite{O}, p.~120) induces the residue map
\begin{equation} \label{residueX}
\Omega_{X^{\dagger}/k^{\dagger}} \map \oplus_i \, \cO_{D_i}. 
\end{equation}
With notation as in \ref{points},
the $(k-1)$th exterior power of the residue map (\ref{residueX}) is the iterated residue map
\begin{equation}\label{iteratedresidueX}
\omega_X(D) \map \oplus_{i_1 < \cdots < i_{k-1}} \, k(P_{i_1,\ldots,i_{k-1}}) 
\end{equation}
given locally at $P=P_{i_1,\ldots,i_{k-1}}$ by
$$ \frac{dx_1}{x_1} \wedge \cdots \wedge \frac{dx_{k-1}}{x_{k-1}} \mapsto 1.$$
Taking global sections, we obtain a natural map
$$H^0(\omega_X(D)) \map \wedge^{k-1}\bC^n.$$

\begin{thm} \label{embcan}
The embedding $(X,D) \inj G(k-1,n-1)_e$ is canonical. Specifically,
the composition $$(X,D) \inj G(k-1,n-1)_e \inj \bP(\wedge^{k-1} \mathfrak{h})$$
is defined by the line bundle $\omega_X(D)$ together with a natural isomorphism
$$H^0(\omega_X(D)) \stackrel{\sim}{\map} \wedge^{k-1}\mathfrak{h}^{\vee} \subset \wedge^{k-1} \bC^n$$
given by the iterated residue map (\ref{iteratedresidueX}).
\end{thm}

\begin{proof}
The map $\mathfrak{h}^{\vee} \map H^0(\Omega_{X^{\dagger}/k^{\dagger}})$ is given by 
$m \mapsto d\chi^m/\chi^m|_X$ for $m \in M \subset \mathfrak{h}^{\vee}$. 
The composition
$$\mathfrak{h}^{\vee} \map H^0(\Omega_{X^{\dagger}/k^{\dagger}}) \map \oplus_i \, H^0(\cO_{D_i}) = {\bC^n}$$
is the embedding $\mathfrak{h}^{\vee}=(\sum x_i = 0) \inj \bC^n$.
For, the divisor $B_i^+ \subset Z$ corresponds to the facet $\Gamma_i^+$ of $\Delta(k,n) \subset M_{\bR}$, 
with outward normal $e_i \in N_{\bR}$. Hence the log differential $d\chi^m/\chi^m$ on $Z$ has residue 
$\langle m,e_i \rangle$
along $B_i^+$, the order of vanishing of the character $\chi^m$ along $B_i^+ \subset Z$ (\cite{F}, p.~61). 
Consequently, the composition
$$\wedge^{k-1}\mathfrak{h}^{\vee} \map H^0(\omega_X(D)) \map \wedge^{k-1}\bC^n$$
is the embedding $\wedge^{k-1}\mathfrak{h}^{\vee} \inj \wedge^{k-1}\bC^n$.
We prove below (Proposition~\ref{canbasis}) that the map 
$H^0(\omega_X(D)) \map \wedge^{k-1}\bC^n$ is an isomorphism onto 
$\wedge^{k-1}~\mathfrak{h}^{\vee} \inj \bC^n$. 
Hence the map $\wedge^{k-1} \mathfrak{h}^{\vee} \map H^0(\omega_X(D))$
yielding the embedding $(X,D) \inj G_e \inj \bP(\wedge^{k-1} \mathfrak{h})$ is given by the inverse of this 
map, and the embedding is uniquely determined as desired. 
\end{proof}

\begin{lem} \label{adjunction}
For $k \ge 3$, the natural map $H^0(\omega_X(D)) \map H^0(\omega_D)$ is an isomorphism.
\end{lem}
\begin{proof}
By adjunction, we have an exact sequence
$$0 \map \omega_X \map \omega_X(D) \map \omega_D \map 0.$$
Consider the associated long exact sequence of cohomology
$$0 \map H^0(\omega_X) \map H^0(\omega_X(D)) \map H^0(\omega_D) \map H^1(\omega_X) \map \cdots.$$
We have $H^0(\omega_X)=H^{k-1}(\cO_X)^{\vee}$ and $H^1(\omega_X)=H^{k-2}(\cO_X)^{\vee}$ by Serre duality,
and $H^{k-1}(\cO_X)=H^{k-2}(\cO_X)=0$ by Proposition~\ref{rationality}, hence the map
$H^0(\omega_X(D)) \map H^0(\omega_D)$ is an isomorphism as claimed.
\end{proof}

\begin{lem} \label{inductionexseq}
Write $D_{ij} = D_i \cap D_j$.
For $k \ge 3$, the sequence
$$0 \map H^0(\omega_X(D)) \map \oplus_i \, H^0(\omega_{D_i}(D-D_i)) \map \oplus_{i<j} \, H^0(\omega_{D_{ij}}(D-D_i-D_j))$$
is exact.
\end{lem}
\begin{proof}
Consider the exact sequence
$$0 \map \cO_D \map \oplus_i \, \cO_{D_i} \map \oplus_{i<j} \, \cO_{D_{ij}}$$
of sheaves on $D$ given by the decomposition $D=\sum D_i$.
Tensoring with the line bundle $\omega_X(D)$, we obtain the exact sequence
$$0 \map \omega_D \map  \oplus_i \, \omega_{D_i}(D-D_i) \map \oplus_{i<j} \, \omega_{D_{ij}}(D-D_i-D_j)$$
by adjunction. Taking global sections and applying Lemma~\ref{adjunction}, we obtain our result.
\end{proof}

\begin{prop} \label{canbasis}
The map $H^0(\omega_X(D)) \map \wedge^{k-1} \bC^n$ is an isomorphism onto 
$\wedge^{k-1} \mathfrak{h}^{\vee} \inj \wedge^{k-1} \bC^n$.
\end{prop}
\begin{proof}
The proof is by induction on $k$.
The case $k=1$ is trivial ($X$ is a point) and the case $k=2$ is known (\cite{K2}, p.~242).
So, assume that $k \ge 3$ and the result holds for $k-2$ and $k-1$.
Consider the exact sequence 
$$0 \map H^0(\omega_X(D)) \map \oplus_i \, H^0(\omega_{D_i}(D-D_i)) \stackrel{\theta}{\map} \oplus_{i<j} \, H^0(\omega_{D_{ij}}(D-D_i-D_j))$$
of Lemma~\ref{inductionexseq}.
By Proposition~\ref{restriction}, the pair $(D_i,(D-D_i)|_{D_i})$ is a fibre of the family over the Chow quotient 
$G(k-1,n-1)_{e_i} \cq  H_i$, where $H_i$ is the quotient of $H$ by the $i$th coordinate $\bC^{\times}$.
Similiarly, $(D_{ij}, (D-D_i-D_j)|_{D_{ij}})$ is a fibre of the family over $G(k-2,n-2)_{e_i,e_j} \cq H_{ij}$.
Write $\mathfrak{h}_i= \Lie H_i = \bC^n / \langle e, e_i \rangle$
and $\mathfrak{h}_{ij}= \Lie H_{ij} =\bC^n / \langle e,e_i,e_j \rangle$.
By induction, the map $\theta$ is identified with the map 
\begin{equation} \label{inductionmap}
\oplus_i \, \wedge^{k-2} \mathfrak{h}_i^{\vee} \map \oplus_{i<j} \, \wedge^{k-3} \mathfrak{h}_{ij}^{\vee}
\end{equation}
given by
$$\wedge^{k-2}\mathfrak{h}_i^{\vee} \map \wedge^{k-3} \mathfrak{h}_{ij}^{\vee}, \ e_j^{\star} \wedge \omega \mapsto \omega.$$
Here $e_1^{\star}, \cdots, e_n^{\star}$ denotes the basis of ${\bC^n}^{\vee}$ dual to the standard basis $e_1,\ldots,e_n$ of $\bC^n$. 
Hence the kernel $H^0(\omega_X(D))$ of $\theta$ is identified with the kernel of the map (\ref{inductionmap}).
To compute this kernel, we may replace $(X,D)$ by a generic pair $(X',D')$, i.e., a hyperplane arrangement.
Then $H^0(\omega_{X'}(D'))$ is identified with $\wedge^{k-1}\mathfrak{h}^{\vee}$ via the residue map 
(cf. \cite{K1}, p.~71, Proposition~3.2.3(a)), equivalently, the kernel of (\ref{inductionmap}) is the image of the embedding
$$\wedge^{k-1} \mathfrak{h}^{\vee} \inj \oplus_i \, \wedge^{k-2}\mathfrak{h}_i^{\vee}$$
given by
$$\wedge^{k-1} \mathfrak{h}^{\vee} \map \wedge^{k-2}\mathfrak{h}_i^{\vee}, \ e_i^{\star} \wedge \omega \mapsto \omega.$$ 
This completes the proof.
\end{proof}

\begin{lem} \label{Zstratadecomp}
Let $Z$ be the degenerate toric variety corresponding to a polyhedral subdivision $\underline{Q}$ of a lattice polytope $Q$. 
Let $Z^j$ denote the disjoint union of the toric strata of $Z$ of codimension $j$ which are not contained
in the toric boundary $B$; write $p^j \colon Z^j \map Z$ for the induced map.
There is an exact sequence
$$0 \map \cO_Z \map p^0_{\star} \cO_{Z^0} \map p^1_{\star} \cO_{Z^1} \map \cdots.$$
\end{lem}

\begin{proof}
We define a complex 
$$0 \map \cO_{Z} \map p^0_{\star} \cO_{Z^0} \map p^1_{\star} \cO_{Z^1} \map \cdots $$
as follows. The map $p^i_{\star} \cO_{Z^i} \map p^{i+1}_{\star} \cO_{Z^{i+1}}$ is given by the restriction 
maps from codimension $i$ strata to codimension $i+1$ strata with various signs.
The signs are determined by fixing an orientation for each polytope $\sigma \in \underline{Q}$.
We assume that each maximal polytope has orientation induced by a fixed orientation of $Q$,
then the map $\cO_Z \map p^0_{\star} \cO_{Z^0}$ is given by restriction (no signs). 
Let $Z(\sigma)$ denote the stratum of $Z$ corresponding to a face $\sigma$ of $\underline{Q}$.
For $\sigma$ a face of $\underline{Q}$ and $\tau$ a facet of $\sigma$ 
the map $\cO_{Z(\sigma)} \map \cO_{Z(\tau)}$ has sign $+1$ if $\sigma$ and $\tau$ are oriented compatibly
and $-1$ otherwise. 
 
We show that the corresponding complex of graded rings 
\begin{equation} \label{seq-R}
0 \map R_Z \map R_{Z^0} \map R_{Z^1} \map \cdots
\end{equation}
is exact. Recall that, as a $\bC$-vector space,
$$R_Z=\langle \chi^s \ | \ s \in \Cone(Q) \cap M' \rangle_{\bC}$$
where $M' = \bZ \oplus M$, and $Q$ sits inside the affine hyperplane $(1,M_{\bR})$ of $M'_{\bR}$.
Similiarly,
$$R_{Z^j}=\bigoplus_{Q' \in \underline{Q}(j)} \langle \chi^s \ | \ s \in \Cone(Q') \cap M' \rangle_{\bC}$$
where $\underline{Q}(j)$ is the set of codimension $j$ faces of $\underline{Q}$ which are not contained in the boundary.
The complex (\ref{seq-R}) is a direct sum of complexes, with one summand for each $s \in \Cone(Q) \cap M'$.
The summand corresponding to $s$ is
\begin{equation} \label{seq-s} 
0 \map \bC \map \bigoplus_{Q' \in \underline{Q}(0,s)} \bC \map \bigoplus_{Q' \in \underline{Q}(1,s)} \bC \map \cdots
\end{equation}
where $\underline{Q}(j,s)=\{Q' \in \underline{Q}(j) \ | \ s \in \Cone(Q') \}$.
We check below that each such summand is exact, completing the proof of (1).

Consider the complex
\begin{equation} \label{seq-s'}
0 \map \bigoplus_{Q' \in \underline{Q}(0,s)} \bC \map \bigoplus_{Q' \in \underline{Q}(1,s)} \bC \map \cdots
\end{equation}
obtained from (\ref{seq-s}) by omitting the first term.
Assume $s \neq 0$, let $v$ be the point of $Q$ given by the ray $\bR_{\ge 0} \cdot s \subset \Cone(Q)$
and write $(\underline{Q}-v)$ for the subcomplex of $\underline{Q}$ consisting of cells not containing $v$.
The complex (\ref{seq-s'}) is identified with the chain complex 
$C_{\cdot}(\underline{Q},(\underline{Q}-v) \cup \underline{\partial Q})$
computing the cellular homology of the pair of CW-complexes $(\underline{Q},(\underline{Q}-v) \cup \underline{\partial Q})$.
But the pair $(\underline{Q},(\underline{Q}-v) \cup \underline{\partial Q})$ is homotopy equivalent to 
the pair $(B^d,S^{d-1})$ consisting of the $d$-dimensional ball and its boundary; here $d= \dim Q$.
It follows that the complex (\ref{seq-s}) is exact as required.
Similiarly, if $s=0$, the complex (\ref{seq-s'}) is identified with $C_{\cdot}(\underline{Q},\underline{\partial Q})$
and so (\ref{seq-s}) is again exact.
\end{proof}

\begin{lem} \label{Xstratadecomp}
Let $(X,D)$ be the stable pair associated to $[Z] \in G(k,n)\cq H$.
Consider the stratification of $X$ induced by the stratification of $Z$ by orbit closures.
Let $X^j$ denote the disjoint union of the strata of $X$ of codimension $j$ which are not contained
in $D$; write $p^j \colon X^j \map X$ for the induced map.
There is an exact sequence
$$0 \map \cO_X \map p^0_{\star} \cO_{X^0} \map p^1_{\star} \cO_{X^1} \map \cdots.$$
\end{lem}

\begin{proof}
The result follows immediately from Lemma~\ref{Zstratadecomp} and the smoothness of the map 
$f_Z \colon H \times X \map Z$ (Proposition~\ref{(X,D)slc}).
\end{proof}

\begin{prop} \label{rationality}
Let $(X,D)$ be the stable pair associated to $[Z] \in G(k,n)\cq H$.
Then $H^i(\cO_X)=0$ for $i > 0$.
\end{prop}
\begin{proof}
We use the exact sequence
$$0 \map \cO_X \map p^0_{\star} \cO_{X^0} \map p^1_{\star} \cO_{X^1} \map \cdots$$
of Lemma~\ref{Xstratadecomp}. Each stratum of $X$ is rational and has rational singularities by
Theorem~\ref{Xstrata}. Hence $H^i(\cO_{X^j})=0$ for each $i > 0$, $j \ge 0$.
So, the group $H^i(\cO_X)$ is the $i$th cohomology group of the complex
\begin{equation} \label{globalseq}
0 \map H^0(\cO_{X^0}) \map H^0(\cO_{X^1}) \map \cdots.
\end{equation}
Let $\underline{Q}$ denote the polyhedral subdivision of $Q=\Delta(k,n)$ given by $Z$.
Recall that each stratum of $Z$ which is not contained in the toric boundary yields a stratum of $X$ of the same codimension.
Thus, the complex (\ref{globalseq}) is identified with the chain complex $C_{\cdot}(\underline{Q}, \underline{\partial Q})$
computing the cellular homology of the pair of CW-complexes $(\underline{Q}, \underline{\partial Q})$ 
(cf. proof of Lemma~\ref{Zstratadecomp}). Hence $H^i(\cO_X)=0$ for $i>0$ as required.
\end{proof}
  
\section{Proof of the Main Theorem}

\begin{proof}[Proof of Theorem~\ref{mainthm}]
The fibre of the family $(\cX^u,\sD^u)/C_{k,n}$ over a point of $H_{k,n}$ is identified with the 
corresponding hyperplane arrangement by Proposition~\ref{genericfibres}.
For an arbitrary fibre $(X,D)$, the pair $(X,D)$ is slc by Proposition~\ref{(X,D)slc} and 
the sheaf $\omega_X(D)$ is a very ample line bundle by Theorem~\ref{logdiffX}. 
Hence $(X,D)$ is a stable pair. 

We show that no two fibres of $(\cX^u,\sD^u)/C_{k,n}$ are isomorphic.
For each fibre $(X,D)$, the embedding $(X,D) \inj G(k-1,n-1)_e$ is canonical
by Theorem~\ref{embcan}. If $(X,D) \inj G_e$ is the fibre over $[Z] \in G(k,n)\cq H$,
then $Z=\overline{HX} \subset G(k,n)$. Hence $(X,D)$ determines $Z \subset G(k,n)$,
as required.

Suppose given  a stable pair $(X,D)$ and  a smoothing $(\cX,\sD)/(0 \in T)$ of $(X,D)$  
as in Definition~\ref{def_stable_pair}. 
We show that $(X,D)$ occurs as
a fibre of $(\cX^u,\sD^u)/C_{k,n}$. Write $T^{\times} = T \backslash \{ 0 \}$ and let
$(\cX^{\times},\sD^{\times})/T^{\times}$ denote the restriction of $(\cX,\sD)/T$.
Then there is a canonical embedding 
$$(\cX^{\times},\sD^{\times}) \inj G(k-1,n-1)_e \times T^{\times}$$ 
defined by the logarithmic differentials, and $(\cX^{\times},\sD^{\times})/T^{\times}$ is obtained as the pullback 
of the family $(\cX^u,\sD^u)/C_{k,n}$ by some map $T^{\times} \map C_{k,n}$ (\cite{K1}, p.~74, Corollary~3.3.11). 
This map extends uniquely to a map $T \map C_{k,n}$;
let $(\cX',\sD')/T$ denote the family obtained by pullback, and $(X',D')$ the special fibre.
Then, by separatedness of moduli of stable pairs (\cite{A2}, p.~7, Theorem~2.9), $(X,D)$ is isomorphic to $(X',D')$.
Here we require the following `inversion of adjunction' result: $(X,D)$ slc implies $(\cX,X+\sD)$ log canonical.
It is a consequence of the minimal model program in dimension $k=\dim X + 1$ (\cite{A2}, p.~6, 2.7).
\end{proof}

\begin{rem} \label{MMP}
We note that the minimal model program assumption is only used to guarantee an
inversion of adjunction result, which it may be possible to prove by other methods. 
For example, in some cases, it can be proved using vanishing (\cite{KM}, p.~174, Theorem~5.50).
\end{rem}

\end{document}